\documentclass[12pt, reqno]{amsart}
\usepackage{epsfig}
\usepackage{enumitem}
\usepackage{multicol}
\def\footnoterule{\kern-3pt \hrule width 1.9 cm \kern2.6pt}
\begin{document}
\centerline{\bf \small{ON HOLOMORPHY OF FENYVES BCI-ALGEBRAS}}
\vskip 0.4 cm
\centerline{\scriptsize{E. ILOJIDE\footnote{Corresponding
author}, T. G. JAIY\'E\d OL\'A AND M. O. OLATINWO}}
\vskip 0.4 cm
\begin{center}
\begin{minipage}[b]{9.0cm}
\scriptsize{
\noindent 
ABSTRACT. Fenyves BCI-algebras are BCI-algebras that satisfy the Bol-Moufang identities. In this paper, the holomorphy of BCI-algebras are studied. It is shown that whenever a loop and its holomorph are BCI-algebras, the former is $p$-semisimple if and only if the latter is $p$-semisimple. Whenever a loop and its holomorph are BCI-algebras, it is established that the former is a BCK-algebra if and only if the latter has a BCK-subalgebra. Moreover, the holomorphy of the associative and some non-associative Fenyves BCI-algebras are also studied.
}
\end{minipage}
\end{center}
\par
\vskip 0.5 cm
{\small{
\centerline{1. INTRODUCTION}}
}
\vskip 0.3 cm
\noindent 
BCK-algebras and BCI-algebras are abreviated to two B-algebras. They were introduced by Imai and Iseki \cite{47}.
The two algebras are originated from two different sources. One of the
motivations is based on set theory. In set theory, there are three
most elementary and fundamental operations. They are the union,
intersection, and set difference. If we consider those three
operations and their properties, then as a generalization of them,
we have the notion of Boolean algebras. If we take both the union
and intersection, then as a general algebra, the notion of
distributive lattices is obtained. Moreover, if we consider the
union or the intersection alone, we have the notion of upper
semilattices or lower semilattices. However, the set difference
together with its properties had not been considered systematically
before Imai and Iseki.
Jaiy\'e\d ol\'a et al. \cite{FBCI1} introduced new kinds of BCI-algebras known as Fenyves BCI-algebras. In this paper, we will study the holomorphy of the classical BCI-algebras as well as the holomorphy of Fenyves BCI-algebras.

\vskip 0.5 cm
{\small{
\centerline{2. PRELIMINARY}}
}
\vskip 0.3 cm
\noindent

\par
\noindent
{\bf Definition 1(\cite{2}):} A triple $(X;\ast,0)$ is called a BCI-algebra if the following
conditions are satisfied for any $x,y,z\in X$
\begin{enumerate}
  \item $((x\ast y)\ast (x\ast z))\ast (z\ast y)=0$;
  \item $x\ast 0=x$;
  \item $x\ast y=0$ and $y\ast x=0$ $\Longrightarrow$ $x=y$.
\end{enumerate}
\par

\noindent
We call the binary operation $\ast$ on $X$ the multiplication
on $X$, and the constant $0$ of $X$ the zero element of $X$. We
often write $X$ instead of $(X;\ast,0)$ for a BCI-algebra in
brevity.
\par

\noindent
{\bf Example 1(\cite{2}):} Let $S$ be a set. Denote $2^{S}$ for the power set of $S$ in the
sense that $2^{S}$ is the collection of all subsets of $S$, $-$ the set
difference and $\emptyset$ for the empty set. Then, $(2^{S};-,\emptyset)$ is a
BCI-algebra.
\par

\noindent
{\bf Example 2(\cite{2}):} Suppose $(G; \cdot ,e)$ is an abelian group with $e$ as the identity
element. Define a binary operation $\ast$ on $G$ by $x\ast
y=xy^{-1}$. Then $(G,\ast ,e)$ is a BCI-algebra.
\par

\noindent
{\bf Example 3(\cite{2}):} Let $\mathbb{Z}$ be the set of integers. Then, $(\mathbb{Z},-,0)$ is a BCI-algebra.
\par

\noindent
The following theorem gives necessary and sufficient conditions for
the existence of a BCI-algebra.
\par

\noindent
{\bf Theorem 1(\cite{2}):} Let $X$ be a non-empty set,$\ast$ a binary operation on $X$ and $0$
a constant element of $X$. Then, $(X,\ast,0)$ is a BCI- algebra if
and only if the following conditions hold:
\begin{enumerate}
  \item $((x\ast y)\ast (x\ast z))\ast (z\ast y)=0$;
  \item $(x\ast (x\ast y)\ast y=0$;
  \item $x\ast x=0$;
  \item $x\ast y=0$ and $y\ast x=0$ imply $x=y$.
\end{enumerate}
\par

\noindent
\vskip 0.5 cm
{\small{
\centerline{3. THE HEART OF THE MATTER} }}
\vskip 0.3 cm
\noindent
\par			

\noindent				
{\bf Definition 2:}  A BCI- algebra $(X,\ast,0)$ is called a Fenyves BCI-algebra if it
satisfies any of the identities of Bol-Moufang type.
\par

\noindent
The identities of Bol-Moufang type are given below:
\begin{multicols}{2}
\begin{description}
\item[$F_{1}$:] $xy\cdot zx=(xy\cdot z)x$

\item[$F_{2}$:] $xy\cdot zx=(x\cdot yz)x$ (Moufang identity)

\item[$F_{3}$:] $xy\cdot zx=x(y\cdot zx)$

\item[$F_{4}$:] $xy\cdot zx=x(yz\cdot x)$ (Moufang identity)

\item[$F_{5}$:] $(xy\cdot z)x=(x\cdot yz)x$

\item[$F_{6}$:] $(xy\cdot z)x=x(y\cdot zx)$ (extra identity)

\item[$F_{7}$:] $(xy\cdot z)x=x(yz\cdot x)$

\item[$F_{8}$:] $(x\cdot yz)x=x(y\cdot zx)$

\item[$F_{9}$:] $(x\cdot yz)x=x(yz\cdot x)$

\item[$F_{10}$:] $x(y\cdot zx)=x(yz\cdot x)$

\item[$F_{11}$:] $xy\cdot xz=(xy\cdot x)z$

\item[$F_{12}$:] $xy\cdot xz=(x\cdot yx)z$

\item[$F_{13}$:] $xy\cdot xz=x(yx\cdot z)$ (extra identity)

\item[$F_{14}$:] $xy\cdot xz=x(y\cdot xz)$

\item[$F_{15}$:] $(xy\cdot x)z=(x\cdot yx)z$

\item[$F_{16}$:] $(xy\cdot x)z=x(yx\cdot z)$

\item[$F_{17}$:] $(xy\cdot x)z=x(y\cdot xz)$ (Moufang identity)

\item[$F_{18}$:] $(x\cdot yx)z=x(yx\cdot z)$

\item[$F_{19}$:] $(x\cdot yx)z=x(y\cdot xz)$ (left Bol identity)

\item[$F_{20}$:] $x(yx\cdot z)=x(y\cdot xz)$

\item[$F_{21}$:] $yx\cdot zx=(yx\cdot z)x$

\item[$F_{22}$:] $yx\cdot zx=(y\cdot xz)x$ (extra identity)

\item[$F_{23}$:] $yx\cdot zx=y(xz\cdot x)$

\item[$F_{24}$:] $yx\cdot zx=y(x\cdot zx)$

\item[$F_{25}$:] $(yx\cdot z)x=(y\cdot xz)x$

\item[$F_{26}$:] $(yx\cdot z)x=y(xz\cdot x)$ (right Bol identity)

\item[$F_{27}$:] $(yx\cdot z)x=y(x\cdot zx)$ (Moufang identity)

\item[$F_{28}$:] $(y\cdot xz)x=y(xz\cdot x)$

\item[$F_{29}$:] $(y\cdot xz)x=y(x\cdot zx)$

\item[$F_{30}$:] $y(xz\cdot x)=y(x\cdot zx)$

\item[$F_{31}$:] $yx\cdot xz=(yx\cdot x)z$

\item[$F_{32}$:] $yx\cdot xz=(y\cdot xx)z$

\item[$F_{33}$:] $yx\cdot xz=y(xx\cdot z)$

\item[$F_{34}$:] $yx\cdot xz=y(x\cdot xz)$

\item[$F_{35}$:] $(yx\cdot x)z=(y\cdot xx)z$

\item[$F_{36}$:] $(yx\cdot x)z=y(xx\cdot z)$ (RC identity)

\item[$F_{37}$:] $(yx\cdot x)z=y(x\cdot xz)$ (C identity)

\item[$F_{38}$:] $(y\cdot xx)z=y(xx\cdot z)$

\item[$F_{39}$:] $(y\cdot xx)z=y(x\cdot xz)$ (LC identity)

\item[$F_{40}$:] $y(xx\cdot z)=y(x\cdot xz)$

\item[$F_{41}$:] $xx\cdot yz=(x\cdot xy)z$ (LC identity)

\item[$F_{42}$:] $xx\cdot yz=(xx\cdot y)z$

\item[$F_{43}$:] $xx\cdot yz=x(x\cdot yz)$

\item[$F_{44}$:] $xx\cdot yz=x(xy\cdot z)$

\item[$F_{45}$:] $(x\cdot xy)z=(xx\cdot y)z$

\item[$F_{46}$:] $(x\cdot xy)z=x(x\cdot yz)$ (LC identity)

\item[$F_{47}$:] $(x\cdot xy)z=x(xy\cdot z)$

\item[$F_{48}$:] $(xx\cdot y)z=x(x\cdot yz)$ (LC identity)

\item[$F_{49}$:] $(xx\cdot y)z=x(xy\cdot z)$

\item[$F_{50}$:] $x(x\cdot yz)=x(xy\cdot z)$

\item[$F_{51}$:] $yz\cdot xx=(yz\cdot x)x$

\item[$F_{52}$:] $yz\cdot xx=(y\cdot zx)x$

\item[$F_{53}$:] $yz\cdot xx=y(zx\cdot x)$ (RC identity)

\item[$F_{54}$:] $yz\cdot xx=y(z\cdot xx)$

\item[$F_{55}$:] $(yz\cdot x)x=(y\cdot zx)x$

\item[$F_{56}$:] $(yz\cdot x)x=y(zx\cdot x)$ (RC identity)

\item[$F_{57}$:] $(yz\cdot x)x=y(z\cdot xx)$ (RC identity)

\item[$F_{58}$:] $(y\cdot zx)x=y(zx\cdot x)$

\item[$F_{59}$:] $(y\cdot zx)x=y(z\cdot xx)$

\item[$F_{60}$:] $y(zx\cdot x)=y(z\cdot xx)$
\end{description}
\end{multicols}
\par

\noindent
Consequent upon this definition, there are sixty varieties of
Fenyves BCI -algebras. We give some varieties of Fenyves
BCI-algebras as follows:
\par

\noindent
{\bf Definition 3:} A BCI-algebra $(X,\ast,0)$ is called an $F_{1}$-algebra if it
satisfies the following condition: $xy\ast zx=(xy\ast z)\ast x
~\forall~ x,y,z \in X$.
\par

\noindent
{\bf Definition 4:} A BCI-algebra $(X,\ast,0)$ is called an $F_{5}$- algebra if it
satisfies the following condition: $(xy\ast z)x=(x\ast yz)x~\forall~
x,y,z \in X$.
\par

\noindent
{\bf Definition 5:} A BCI-algebra $(X,\ast,0)$ is called an $F_{47}$-algebra if it
satisfies the following condition: $(x\ast xy)z=x(xy\ast z)~\forall~
x,y,z\in X$.
\par

\noindent
Let us now give some examples of a Fenyves' BCI-algebra:
\par

\noindent
{\bf Example 4:} Let $(G; \cdot ,e)$ be an abelian group with $e$ as the identity
element. Define a binary operation $\ast$ on $G$ by $x\ast
y=xy^{-1}$. Then $(G,\ast ,e)$ is an $F_{8}$-algebra, $F_{19}$-algebra, $F_{29}$-algebra, $F_{39}$-algebra, $F_{46}$-algebra, $F_{52}$-algebra, $F_{54}$-algebra, $F_{59}$-algebra.
\par

\noindent
{\bf Example 5:} Let $S$ be a set. Denote $2^{S}$ for the power set of $S$ in the
sense that $2^{S}$ is the collection of all subsets of $S$, $-$ the set
difference and $\emptyset$ for the empty set. Then $(2^{S};-,\emptyset)$ is an $F_{5}$-algebra, $F_{42}$-algebra, $F_{54}$-algebra.
\par

\noindent
{\bf Definition 6(\cite{2}):} A BCI-algebra $(X,\ast,0)$ is called associative if $(x\ast y)\ast
z=x\ast (y\ast z)~\forall~ x,y,z\in X$.
\par

\noindent
{\bf Definition 7(\cite{2}):} A BCI-algebra $(X,\ast,0)$ is called $p$-semisimple if  $0\ast (0\ast   x) = x~\forall~ x\in X$ .
\par

\noindent
The following theorems give equivalent conditions for associativity and $p$-semisimplicity in a BCI-algebra:
\par

\noindent
{\bf Theorem 2 \label{conclude3} (\cite{2}):} Given a BCI-algebra $X$, the following are equivalent:
\begin{enumerate}
\item $X$ is associative.
\item $0\ast x=x ~\forall~ x\in X$.
\item $x\ast y=y\ast x ~\forall~ x,y \in X$.
\end{enumerate}
\par

\noindent
{\bf Theorem 3\label{conclude1} (\cite{2}):} Let $X$ be a BCI-algebra. Then the following conditions are
equivalent for any $x,y,z,u\in X$:
\begin{enumerate}
\item $X$ is $p$-semisimple.
\item $(x\ast y)\ast (z\ast u)=(x\ast z)\ast (y\ast u)$.
\item $0\ast (y\ast x)=x\ast y$.
\item $(x\ast y)\ast (x\ast z)=z\ast y$.
\item $z\ast x=z\ast y$ implies $x=y$. (the left cancellation law)
\item $x\ast y=0$ implies $x=y$.
\end{enumerate}
\par

\noindent
{\bf Theorem 4\label{conclude2} (\cite{2}):} Let $X$ be a BCI-algebra. $X$ is $p$-semisimple if and only if one
of the following conditions holds for any $x,y,z\in X$:
\begin{enumerate}
\item $x\ast z=y\ast z$ implies $x=y$. (the right cancellation law)
\item $(y\ast x)\ast (z\ast x)=y\ast z$.
\item $(x\ast y)\ast (x\ast z)=0\ast (y\ast z)$.
\end{enumerate}
\par

\noindent
{\bf Definition 8 (Jaiy\'e\d ol\'a \cite{davidref:10}):} Let $L$ be a non-empty set. Define a binary operation ($\cdot $) on
$L$ . If $x\cdot y\in L$ for all $x, y\in L$, $(L, \dot )$ is
called a groupoid. If the equations:
\begin{displaymath}
a\cdot x=b\qquad\textrm{and}\qquad y\cdot a=b
\end{displaymath}
have unique solutions for $x$ and $y$ respectively, then $(L, \cdot
)$ is called a quasigroup. If there exists a unique element $e\in L$ called the identity element
such that for all $x\in L$, $x\cdot e=e\cdot x=x$, $(L, \cdot )$ is
called a loop.
\par

\noindent
{\bf Theorem 5\label{conclude2.1} (\cite{2}):} Suppose that $(X;\ast , 0)$ is a BCI-algebra. $X$ is associative
if and only if $X$ is $p$-semisimple and $X$ is quasi-associative.
\par

\noindent
{\bf Theorem 6\label{conclude2.2} (\cite{2}):} Suppose that $(X;\ast , 0)$ is a BCI-algebra. Then $(x\ast y)\ast z=(x\ast z)\ast y$ for all $x,y,z\in X$.
\par

\noindent
{\bf Remark 1:} In Theorem 5, quasi-associativity in BCI-algebra plays a similar role which weak associativity (i.e. the $F_i$ identities) plays in quasigroup and loop theory.
\par

\noindent
{\bf Definition 9(\cite{2}):} A BCI-algebra $(X;\ast , 0)$ is called a BCK-algebra if $0\ast x=0$ for all $x\in X$.
\par

\paragraph{}
Since late 1970s, BCI and BCK algebras have been paid much attention. In
particular, the participation in the research of polish
mathematicians Tadeusz Traczyk and Andrzej Wronski as well as
Australian mathematician William H. Cornish, etc. has made
the study of BCI-algebras to attract interest of many mathematicians. Many interesting and
important results are discovered continuously. Now, the theory of
BCI-algebras has been widely spread to areas such as general theory
which include congruences, quotient algebras, BCI-Homomorphisms,
direct sums and direct products, commutative BCK-algebras, positive implicative and implicative BCK-algebras, derivations of BCI-algebras,
and ideal theory of BCI-algebras (\cite{47,28,new1000,new1,new500}).

\paragraph{}
Jaiy\'e\d ol\'a \emph{et al} \cite{FBCI1} looked at Fenyves identities on the
platform of BCI-algebras. They classified the Fenyves BCI-algebras into $46$ associative and
$14$ non-associative types and showed  that some Fenyves identities played the role of quasi-associativity, vis-a-vis Theorem~\ref{conclude2.1} in BCI-algebras. They clarified the relationship between a BCI-algebra, a quasigroup and a loop.
Some of their results are stated below.

\noindent
{\bf Theorem 7\label{bci:group} (\cite{FBCI1}):} \begin{enumerate}
  \item A BCI algebra $X$ is a quasigroup if and only if it is $p$-semisimple.
  \item A BCI algebra $X$ is a loop if and only if it is associative.
  \item An associative BCI algebra $X$ is a Boolean group.
  \end{enumerate}
\par

\noindent
{\bf Theorem 8\label{bci:group2} \cite{FBCI1}:} Let $(X;\ast,0)$ be a BCI-algebra. If $X$ is any of the following Fenyves BCI-algebras, then $X$ is associative.
\begin{multicols}{3}
\begin{enumerate}
\item $F_{1}$-algebra
\item $F_{2}$-algebra
\item $F_{4}$-algebra
\item $F_{6}$-algebra
\item $F_{7}$-algebra
\item $F_{9}$-algebra
\item $F_{10}$-algebra
\item $F_{11}$-algebra
\item $F_{12}$-algebra
\item $F_{13}$-algebra
\item $F_{14}$-algebra
\item $F_{15}$-algebra
\item $F_{16}$-algebra
\item $F_{17}$-algebra
\item $F_{18}$-algebra
\item $F_{20}$-algebra
\item $F_{22}$-algebra
\item $F_{23}$-algebra
\item $F_{24}$-algebra
\item $F_{25}$-algebra
\item $F_{26}$-algebra
\item $F_{27}$-algebra
\item $F_{28}$-algebra
\item $F_{30}$-algebra
\item $F_{31}$-algebra
\item $F_{32}$-algebra
\item $F_{33}$-algebra
\item $F_{34}$-algebra
\item $F_{35}$-algebra
\item $F_{36}$-algebra
\item $F_{37}$-algebra
\item $F_{38}$-algebra
\item $F_{40}$-algebra
\item $F_{41}$-algebra
\item $F_{43}$-algebra
\item $F_{44}$-algebra
\item $F_{45}$-algebra
\item $F_{47}$-algebra
\item $F_{48}$-algebra
\item $F_{49}$-algebra
\item $F_{50}$-algebra
\item $F_{51}$-algebra
\item $F_{53}$-algebra
\item $F_{57}$-algebra
\item $F_{58}$-algebra
\item $F_{60}$-algebra
\end{enumerate}
\end{multicols}
\par

\noindent
{\bf Remark 2:} All other $F_i$'s which are not mentioned in Theorem 8 were found to be non-associative. Every BCI-algebra is naturally an $F_{54}$ BCI-algebra.
\par

\noindent
{\bf Definition 10 (Jaiy\'e\d ol\'a \cite{davidref:10}):}(Holomorph). Let $(Q,\cdot)$ be a groupoid (quasigroup, loop) and $A(Q)\le
AUM(Q,\cdot)$ be a boolean group of automorphisms of the groupoid
(quasigroup, loop) $(Q,\cdot)$. Let $H=A(Q)\times Q$. Define $\circ$
on $H$ as
\begin{displaymath}
(\alpha,x)\circ(\beta,y)=(\alpha\beta,x\beta\cdot y)~\textrm{for
all}~(\alpha,x),(\beta,y)\in H.
\end{displaymath}
$(H,\circ)$ is a groupoid (quasigroup, loop) and is called the
A-holomorph of $(Q,\cdot)$.
\par

\noindent
{\bf Definition 11 (Jaiy\'e\d ol\'a \cite{davidref:10}):}($\lambda$-regular and $\rho$-regular bijections). A bijection $U$ of a groupoid $(G,\cdot)$ is called $\lambda$-regular if there exists an autotopism $(U,I,U)$ of $(G,\cdot)$ such that $xU\cdot y=(x\cdot y)U~\forall ~x,y\in G$.

A bijection $U$ of a groupoid $(G,\cdot)$ is called $\rho$-regular if there exists an autotopism $(I,U,U)$ of $(G,\cdot)$ such that $x\cdot yU=(x\cdot y)U~\forall~x,y\in G$.
\par

\noindent
\paragraph{}
The holomorph of a loop is a loop according to Bruck \cite{phd15}. Since then, the concept of holomorphy of loops
has caught the attention of some researchers. Interestingly, Ad\'en\'iran \cite{phd11} and Robinson \cite{phd40},
Chein and Robinson\cite{phd17}, Ad\'en\'iran \emph{et al} \cite{phd36}, Chiboka and Solarin \cite{phd19},
\cite{phd20}, Bruck \cite{phd15}, Bruck and Paige \cite{phd16}, Robinson
\cite{40}, Huthnance \cite{phd24}, Ad\'en\'iran et tal \cite{ho1} and, Jaiy\'e\d ol\'a and Popoola \cite{phd41} have
respectively studied the holomorphic structures of Bol/Bruck loops, Moufang loops, central loops,
conjugacy closed loops, inverse property loops, A-loops, extra
loops, weak inverse property loops and generalized Bol loops. Isere \emph{et al} \cite{isere1} looked at the holomorphic structure of Osborn loops.

\paragraph{}
In this present work, the holomorphy of BCI-algebras are studied. The holomorphy of both the associative and non-associative Fenyves BCI-algebras are also studied.
\par

\noindent
\vskip 0.5 cm
{\small{
\centerline{4. MAIN RESULTS} }}
\vskip 0.3 cm
\noindent
\par

\noindent
\vskip 0.5 cm
{\small{
\centerline{4.1. REGULAR BIJECTIONS AND AUTOMORPHISMS OF BCI-ALGEBRAS} }}
\vskip 0.3 cm
\noindent
\par

\noindent
\paragraph{}
We first present a result on regular bijections and automorphisms of BCI-algebras.
\par

\noindent
{\bf Lemma 1:} Let $(G, \ast)$ be a BCI-algebra with $\delta \in $ SYM$(G, \ast)$. Then the following hold:
\begin{enumerate}
\item [(i)]$\delta$ is $\lambda$-regular $\Leftrightarrow \delta R_{x}=R_{x}\delta\Leftrightarrow L_{x\delta}=L_{x}\delta$ for all $x\in G$.

\item [(ii)]$\delta$ is $\rho$-regular $\Leftrightarrow \delta L_{x}=L_{x}\delta \Leftrightarrow R_{x\delta}=R_{x}\delta$ for all $x\in G$.

\item [(iii)]$\delta$ is $\mu$-regular $\Leftrightarrow \delta R_{x}=R_{x\delta}\Leftrightarrow L_{x\delta}=\delta L_{x}$ for all $x\in G$.
\end{enumerate}
\par

\noindent
{\bf Proof:}
\begin{enumerate}
\item[(i)] $\delta$ is $\lambda$-regular $\Leftrightarrow (\delta, I, \delta)\in $
AUT $(G, \ast)\Leftrightarrow y\delta \ast xI=(y\ast x)\delta\Leftrightarrow y\delta R_{x}=yR_{x}\delta \Leftrightarrow \delta R_{x}=R_{x}\delta \Leftrightarrow y\delta R_{x}=yR_{x}\delta \Leftrightarrow y\delta \ast x=(y\ast x)\delta \Leftrightarrow xL_{y\delta}=xL_{y}\delta \Leftrightarrow L_{y\delta}=L_{y}\delta$.

\item [(ii)] $\delta$ is $\rho$-regular $\Leftrightarrow (I, \delta, \delta)\in $ AUT $(G, \ast)\Leftrightarrow xI\ast y\delta=(x\ast y)\delta \Leftrightarrow y\delta L_{x}=yL_{x}\delta \Leftrightarrow \delta L_{x}=L_{x}\delta\Leftrightarrow y \delta L_{x}=yL_{x}\delta\Leftrightarrow x\ast y\delta =(x\ast y)\delta\Leftrightarrow xR_{y\delta}=xR_{y}\delta \Leftrightarrow R_{y\delta}=R_{y}\delta$.

\item [(iii)] $\delta$ is $\mu$-regular with adjoint $\delta'=\delta \Leftrightarrow (\delta, \delta'^{-1}, I)\in $ AUT $(G, \ast)\Leftrightarrow x\delta \ast y\delta'^{-1}=(x\ast y)I\Leftrightarrow x\delta \ast y\delta \delta'^{-1}=x\ast y\delta$ (by replacing $y$ by $y\delta$) $\Leftrightarrow x\delta \ast y=x\ast y\delta \Leftrightarrow x\delta R_{y}=xR_{y\delta}\Leftrightarrow \delta R_{y}=R_{y\delta}\Leftrightarrow x\delta R_{y}=xR_{y\delta}\Leftrightarrow x\delta \ast y=x\ast y\delta\Leftrightarrow yL_{x\delta}=y\delta L_{x}\Leftrightarrow L_{x\delta}=\delta L_{x}$.

\end{enumerate}
\par

\noindent
{\bf Lemma 2:} Let $A$ be an automorphism on a BCI-algebra $(G, \cdot)$. Then the following hold for all $x,y \in G$:
\begin{enumerate}
\item [(i)] $R_{y}A=AR_{yA}$
\item [(ii)] $L_{x}A=AL_{xA}$
\end{enumerate}
\par

\noindent
{\bf Proof:} Suppose $A$ is an automorphism on $(G, \cdot)$. Then $(x\cdot y)A=xA\cdot yA$ for all $x,y \in G$.
\begin{enumerate}
\item [(i)] Fixing $x$, we have $xR_{y}A=xAR_{yA}\Rightarrow R_{y}A=AR_{yA}$.
\item [(ii)] Fixing $y$, we have $yL_{x}A=yAL_{xA}\Rightarrow L_{x}A=AL_{xA}$.
\end{enumerate}
\par

\noindent
\vskip 0.5 cm
{\small{
\centerline{4.2. HOLOMORPHY OF BCI-ALGEBRAS} }}
\vskip 0.3 cm
\noindent
\par

\noindent
We now present results on holomorphy of BCI-algebras.
\par

\noindent
{\bf Theorem 9\label{homo1}:} Let $(G,\cdot,0)$ be a BCI-algebra, and let $(H,\circ)$ be the $A$-holomorph of $(G,\cdot)$. $(H,\circ)$ is a BCI-algebra if and only if $[(x\delta\cdot y\delta)\cdot (x\cdot z\gamma)]\cdot (z\cdot y)=0$ for all $x,y,z\in G, \delta, \gamma\in A(G)$.
\par

\noindent
{\bf Proof:} We shall be using Definition 1. Now, $(I,0)\in H$ Consider $(\alpha,x)\circ (I,0)=(\alpha I,xI\cdot 0)=(\alpha,x)$. So, $(I,0)$ is the zero element of $H$.

Now, let $(\alpha,x),(\beta,y)\in H$. Suppose $(\alpha,x)\circ (\beta,y)=(I,0)$ and $(\beta,y)\circ (\alpha,x)=(I,0)$. Then $(\alpha\beta,x\beta\cdot y)=(I,0)$ and $(\beta\alpha,y\alpha\cdot x)=(I,0)\Leftrightarrow \alpha\beta=I, x\beta\cdot y=0$ and $\beta\alpha=I, y\alpha\cdot x=0\Leftrightarrow \alpha=\beta^{-1}, x\beta\cdot y=0$ and $\beta=\alpha^{-1}, y\alpha\cdot x=0\Leftrightarrow \alpha=\beta, x\beta\cdot y=0$ and $y\alpha\cdot x=0\Leftrightarrow x\alpha\cdot y=0$ and $y\alpha\cdot x=0$. Replace $x$ with $x\alpha^{-1}$ in $x\alpha\cdot y=0$ and $y$ with $y\alpha^{-1}$ in $y\alpha\cdot x=0$ to get $x\cdot y=0$ and $y\cdot x=0\Rightarrow x=y\Rightarrow (\alpha,x)=(\beta,y)$.

Now, let $(\alpha,x),(\beta,y), (\gamma,z)\in H$. Consider $\Big\{\big[(\alpha,x)\circ (\beta,y)\big]\circ\big[(\alpha,x)\circ (\gamma,z)\big]\Big\}\circ \big[(\gamma,z)\circ (\beta,y)\big]=\Big[(\alpha\beta,x\beta\cdot y)\circ (\alpha\gamma,x\gamma\cdot z)\Big]\circ \big(\gamma\beta,z\beta\cdot y\big)=\Big[(\alpha\beta\alpha\gamma,(x\beta\cdot y)\alpha\gamma\cdot (x\gamma\cdot z))\Big]\cdot \big(\gamma\beta,z\beta\cdot y\big)=\Big[\alpha\beta\alpha\gamma\gamma\beta,((x\beta\cdot y)\alpha\gamma\cdot (x\gamma\cdot z))\gamma\beta\cdot (z\beta\cdot y)\Big]$ (as $\alpha,\beta,\gamma$ are homomorphisms)$=[I,((x\beta\alpha\gamma\cdot y\alpha\gamma)\cdot (x\gamma\cdot z)\gamma\beta\cdot (z\beta\cdot y)]=\Big(I,((x\alpha\cdot y\alpha\beta)\cdot (x\beta\cdot z\gamma\beta))\cdot (z\beta\cdot y)\Big)=(I,0)\Leftrightarrow \Big[(x\alpha\cdot y\alpha\beta)\cdot (x\beta\cdot z\beta\gamma)\Big]\cdot (z\beta\cdot y)=0$. Replace $z$ with $z\beta^{-1}$, $x$ with $x\beta^{-1}$ to get $\Big[(x\alpha\beta\cdot y\alpha\beta)\cdot (x\cdot z\gamma)\Big]\cdot (z\cdot y)=0$ (as $\alpha\beta=\beta\alpha$ for Boolean group) $\Leftrightarrow \big[(x\delta \cdot y\delta)\cdot (x\cdot z\gamma)\big]\cdot (z\cdot y)=0$ (where $\delta=\alpha\beta$).\\

Conversely, Suppose $\big[(x\delta\cdot y\delta)\cdot (x\cdot z\gamma)\big]\cdot (z\cdot y)=0$ for all $x,y,z\in G, \delta, \gamma\in A(G)$. Doing the reverse on $\big[(x\alpha\beta\cdot y\alpha\beta)\cdot (x\cdot z\gamma)\big]\cdot (z\cdot y)$ in the equation $\bigg(I,((x\alpha\cdot y\alpha\beta)\cdot (x\beta\cdot z\gamma\beta))\cdot (z\beta\cdot y)\bigg)=(I,0)$ we have $\bigg(\alpha\beta\alpha\gamma\gamma\beta,\Big[(x\alpha\cdot y\alpha\beta)\cdot (x\beta\cdot z\beta\gamma)\Big]\cdot (z\beta\cdot y)\bigg)=(I,0)$ which subsequently leads to $\Big\{[(\alpha,x)\circ (\beta,y)]\circ[(\alpha,x)\circ (\gamma,z)]\Big\}\circ \big[(\gamma,z)\circ (\beta,y)\big]=(I,0)$. Hence, $(H,\circ)$ is a BCI-algebra.
\par

\noindent
{\bf Theorem 10\label{homo2}:} Let $(G,\cdot)$ be a groupoid with $A$-holomorph $(H,\circ)$. $(H,\circ)$ is a BCI-algebra if and only if $(G,\cdot)$ is a BCI-algebra and $[(x\delta\cdot y\delta)\cdot (x\cdot z\gamma)]\cdot (z\cdot y)=0$ for all $x,y,z\in G, \delta, \gamma\in A(G)$.
\par

\noindent
{\bf Proof:} Suppose $(H,\circ)$ is a BCI-algebra. The set $S=\{(I,x):x\in G\}$ is a subalgebra of $(H,\circ)$. Consider the function $f:(G,\cdot)\rightarrow S$ defined by $f(x)=(I,x)\forall x\in G$. We claim that $f$ is an isomorphism.

Now, let $(I,x), (I,y)\in S$. Then $(I,x)\circ (I,y)=(II,xI\cdot y)=(I,x\cdot y)$.

Let $x,y\in G$. Consider $f(x\cdot y)=(I,x\cdot y)=(II,xI\cdot y)=(I,x)\circ (I,y)=f(x)\circ f(y)$. So, $f$ is a homomorphism.

Now, let $x,y\in G\ni f(x)=f(y)\Rightarrow (I,x)=(I,y)\Rightarrow x=y$. So, $f$ is one to one. Clearly, by definition, $f$ is onto. So, $f$ is an isomorphism.

Now, since $(H,\circ)$ is a BCI-algebra, its subalgebra $S$ is also a BCI-algebra, and by the isomorphism, $(G,\cdot)$ is a BCI-algebra.
By Theorem 9, we have $[(x\delta\cdot y\delta)\cdot (x\cdot z\gamma)]\cdot (z\cdot y)=0$.

Conversely, suppose $(G,\cdot)$ is a BCI-algebra and $[(x\delta\cdot y\delta)\cdot (x\cdot z\gamma)]\cdot (z\cdot y)=0$, then by Theorem 9, $(H,\circ)$ is a BCI-algebra.
\par

\noindent
{\bf Theorem 11:} Let $(G,\cdot)$ be a BCI-algebra with an $A$-holomorph $(H,\circ)$ which is also a BCI-algebra. $(H,\circ)$ is $p$-semisimple if and only if $(G,\cdot)$ is $p$-semisimple.
\par

\noindent
{\bf Proof:}
$(H,\circ)$ is $p$-semisimple $\Leftrightarrow (I,0)\circ [(I,0)\circ (\alpha,x)]=(\alpha,x)\Leftrightarrow (I,0)\circ [(I,0)\circ (\alpha,x)]=(\alpha,x)\Leftrightarrow (I,0)\circ [(I\alpha,0\alpha\cdot x)]=(\alpha,x)\Leftrightarrow (I,0)\circ (\alpha,0\cdot x)=(\alpha,x)\Leftrightarrow (I\alpha,0\alpha\cdot 0x)=(\alpha,x)\Leftrightarrow (\alpha,0\cdot (0\cdot x))=(\alpha,x)\Leftrightarrow 0\cdot (0\cdot x)=x\Leftrightarrow (G,\cdot)$ is $p$-semisimple.
\par

\noindent
{\bf Corollary 1:} Let $(G,\cdot)$ be a groupoid with $A$-holomorph $(H,\circ)$. $(H,\circ)$ is a $p$-semisimple BCI-algebra if and only if $(G,\cdot)$ is a $p$-semisimple BCI-algebra and $[(x\delta\cdot y\delta)\cdot (x\cdot z\gamma)]\cdot (z\cdot y)=0$.
\par

\noindent
{\bf Proof:} Follows from Theorem 10 and Theorem 11.
\par

\noindent
{\bf Theorem 12:} \label{homo2.1} Let $(G,\cdot)$ be a BCI-algebra with $A$-holomorph $(H,\circ)$ which is also a BCI-algebra. $(G,\cdot)$ is a BCK-algebra if and only if $\{(I,x):x\in G\}$ is a BCK-algebra.
\par

\noindent
{\bf Proof:} Now, $\{(I,x):x\in G\}$ is BCK $\Leftrightarrow (I,0)\circ (I,x)=(I,0)\Leftrightarrow (II,0I\cdot x)=(I,0)\Leftrightarrow (I,0\cdot x)=(I,0)\Leftrightarrow 0\cdot x=0\Leftrightarrow (G,\cdot)$ is BCK.
\par

\noindent
{\bf Corollary 2:} Let $(G,\cdot)$ be a groupoid with $A$-holomorph $(H,\circ)$. $(H,\circ)$ is a BCK-algebra if and only if $(G,\cdot)$ is a BCK-algebra and $[(x\delta\cdot y\delta)\cdot (x\cdot z\gamma)]\cdot (z\cdot y)=0$ for all $x,y,z\in G, \delta, \gamma\in A(G)$.
\par

\noindent
{\bf Proof:} Follows from Theorem 11 and Theorem 12.
\par

\noindent
\vskip 0.5 cm
{\small{
\centerline{4.3. HOLOMORPHY OF ASSOCIATIVE FENYVES BCI-ALGEBRAS} }}
\vskip 0.3 cm
\noindent
\par

\noindent
We now present results on holomorphy of associative Fenyves BCI-algebras.
\par

\noindent
{\bf Theorem 13:}\label{homo3}
Let $(G,\cdot, 0)$ be a BCI-algebra with an $A$-holomorph $(H,\circ)$ which is also a BCI-algebra. $(H,\circ)$ is associative if and only if $(G,\cdot)$ is associative.
\par

\noindent
{\bf Proof:} Now, $(H,\circ)$ is associative if and only if $(I,0)\circ (\alpha,x)=(\alpha,x)\Leftrightarrow (I\alpha,0\alpha\cdot x)=(\alpha,x)\Leftrightarrow 0\cdot x=x\Leftrightarrow (G,\cdot)$ is associative.
\par

\noindent
{\bf Corollary 3:} Let $(G,\cdot, 0)$ be a BCI-algebra with an $A$-holomorph $(H,\circ)$ which is also a BCI-algebra. $(H,\circ)$ is an $F_i$-algebra if and only if $(G,\cdot)$ is an $F_i$-algebra where $i=1,2, 4, 6, 7, 9, 10, 11, 12, 13, 14,$

$15, 16, 17, 18, 20, 22,
23, 24, 25, 26, 27, 28, 30, 31, 32, 33, 34, 35, 36, 37,$

$38, 40, 41, 43, 44, 45, 47, 48, 49, 50, 51, 53, 57
58, 60$.
\par

\noindent
{\bf Proof:} Follows from Theorem 8 and Theorem 13.
\par

\noindent
{\bf Corollary 4:} Let $(G,\cdot)$ be a groupoid with $A$-holomorph $(H,\circ)$. $(H,\circ)$ is an associative BCI-algebra if and only if $(G,\cdot)$ is an associative BCI-algebra and $[(x\delta\cdot y\delta)\cdot (x\cdot z\gamma)]\cdot (z\cdot y)=0$ for all $x,y,z\in G, \delta, \gamma\in A(G)$.
\par

\noindent
{\bf Proof:} Follows from Theorem 8, Theorem 10 and Theorem 13.
\par

\noindent
\vskip 0.5 cm
{\small{
\centerline{4.4. HOLOMORPHY OF BCI-ALGEBRAS AND NON-ASSOCIATIVE}\centerline{ FENYVES BCI-ALGEBRAS} }}
\vskip 0.3 cm
\noindent
\par

\paragraph{}
We now present results on holomorphy of BCI-algebras and non-associative Fenyves BCI-algebras.

\noindent
{\bf Theorem 14:} Let $(G,\cdot)$ be a BCI-algebra with $A$-holomorph $(H,\circ)$ which is also a BCI-algebra such that every automorphism $\delta$ on $(G,\cdot)$ is $\rho$-regular and $\mid\delta\mid =2$. Then $(H,\circ)$ is an $F_{i}$-algebra if and only if $(G,\cdot)$ is an $F_{i}$-algebra; where $i=27,38$.
\par

\noindent
{\bf Proof:} Let $(H,\circ)$ be an $F_{27}$-algebra, and let $(\alpha,x),(\beta,y),(\gamma,z)\in (H,\circ)$. Then $\{[(\beta,y)\circ (\alpha,x)]\circ (\gamma,z)\}\circ (\alpha,x)=(\beta,y)\circ\{(\alpha,x)\circ [(\gamma,z)\circ (\alpha,x)]\}\Leftrightarrow \{[\beta\alpha,y\alpha\cdot x]\circ (\gamma,z)\}\circ (\alpha,x)=(\beta,y)\circ \{(\alpha,x)\circ [\gamma\alpha,z\alpha\cdot x]\}\Leftrightarrow \{\beta\alpha\gamma,(y\alpha\cdot x)\gamma\cdot z\}\circ (\alpha,x)=(\beta,y)\circ \{\alpha\gamma\alpha,x\gamma\alpha\cdot (z\alpha\cdot x)\}\Leftrightarrow \{\beta\alpha\gamma,(y\alpha\gamma\cdot x\gamma)\cdot z\}\cdot (\alpha,x)=(\beta,y)\circ \{\alpha\gamma\alpha,x\gamma\alpha\cdot (z\alpha\cdot x)\}\Leftrightarrow \{\beta\alpha\gamma\alpha,[(y\alpha\gamma\cdot x\gamma)\cdot z]\alpha\cdot x\}=\{\beta\alpha\gamma\alpha,y\alpha\gamma\alpha\cdot [x\gamma\alpha\cdot (z\alpha\cdot x)]\}\Leftrightarrow \{\beta\alpha\gamma\alpha,[(y\alpha\gamma\cdot x\gamma)\alpha\cdot z\alpha]\cdot x\}=\{\beta\alpha\gamma\alpha,y\alpha\gamma\alpha\cdot [x\gamma\alpha\cdot (z\alpha\cdot x)]\}\Leftrightarrow \{\beta\alpha\gamma\alpha,[(y\alpha\gamma\alpha\cdot x\gamma\alpha)\cdot z\alpha]\cdot x\}=\{\beta\alpha\gamma\alpha,y\alpha\gamma\alpha\cdot [x\gamma\alpha\cdot (z\alpha\cdot x)]\}\Leftrightarrow [(y\alpha\gamma\alpha\cdot x\gamma\alpha)\cdot z\alpha]\cdot x=y\alpha\gamma\alpha\cdot [x\gamma\alpha\cdot (z\alpha\cdot x)]$.

Replace $z$ with $z\alpha^{-1}$, $y$ with $y(\alpha\gamma\alpha)^{-1}$ to get $[(y\cdot x\gamma\alpha)\cdot z]\cdot x=y\cdot [x\gamma\alpha\cdot (z\cdot x)]$.

Put $\delta=\gamma\alpha$ to get $[(y\cdot x\delta)\cdot z]\cdot x=y\cdot [x\delta\cdot (z\cdot x)]$.

Replace $x$ with $x\delta$ to get $[(y\cdot x\delta^{2})\cdot z]\cdot x\delta=y\cdot [x\delta^{2}\cdot (z\cdot x\delta)]\Leftrightarrow [(y\cdot x)\cdot z]\cdot x\delta=y\cdot [x\cdot (z\cdot x\delta)]\Leftrightarrow zL_{yx}R_{x\delta}=zR_{x\delta}L_{x}L_{y}\Leftrightarrow L_{yx}R_{x\delta}=R_{x\delta}L_{x}L_{y}\Leftrightarrow \delta L_{yx}R_{x\delta}=\delta R_{x\delta}L_{x}L_{y}\Leftrightarrow L_{yx}\delta R_{x\delta}=R_{x}\delta L_{x}L_{y}\Leftrightarrow L_{yx}R_{x}\delta=R_{x}L_{x}L_{y}\delta \Leftrightarrow L_{yx}R_{x}=R_{x}L_{x}L_{y}\Leftrightarrow zL_{yx}R_{x}=zR_{x}L_{x}L_{y}$ as required.

Let $(H,\circ)$ be an $F_{38}$-algebra, and let $(\alpha,x),(\beta,y),(\gamma,z)\in (H,\circ)$. Then $[(\beta,y)\circ (I,0)]\circ (\gamma,z)=(\beta,y)\circ [(I,0)\circ (\gamma,z)]\Leftrightarrow [\beta I,yI\cdot 0]\circ (\gamma,z)=(\beta,y)\circ [I\gamma,0\gamma\cdot z]\Leftrightarrow [\beta\gamma,y\gamma\cdot z]=[\beta\gamma,y\gamma\cdot (0\cdot z)]\Leftrightarrow y\gamma\cdot z=y\gamma\cdot (0\cdot z)\Leftrightarrow zL_{y\gamma}=zL_{0}L_{y\gamma}\Leftrightarrow L_{y\gamma}=L_{0}L_{y\gamma}\Leftrightarrow \gamma L_{y\gamma}=\gamma L_{0}L_{y\gamma}\Leftrightarrow L_{y}\gamma =L_{0}\gamma L_{y\gamma}\Leftrightarrow L_{y}\gamma=L_{0}L_{y}\gamma\Leftrightarrow L_{y}=L_{0}L_{y}\Leftrightarrow zL_{y}=zL_{0}L_{y}$ as required.
\par

\noindent
{\bf Theorem 15:} Let $(G,\cdot)$ be a BCI-algebra with $A$-holomorph $(H,\circ)$ which is also a BCI-algebra such that every automorphism $\delta$ on $(G,\cdot)$ is $\lambda$-regular, $\rho$-regular and $\mid\delta\mid =2$. Then $(H,\circ)$ is an $F_{i}$-algebra if and only if $(G,\cdot)$ is an $F_{i}$-algebra; where $i=30,40,50,53,$

$55,56,58$.
\par

\noindent
{\bf Proof:} Let $(H,\circ)$ be an $F_{30}$-algebra, and let $(\alpha,x),(\beta,y),(\gamma,z)\in H$. Then $(\beta,y)\circ \{[(\alpha,x)\circ (\gamma,z)]\circ (\alpha,x)\}=(\beta,y)\circ \{(\alpha,x)\circ [(\gamma,z)\circ (\alpha,x)]\}\Leftrightarrow (\beta,y)\circ \{[\alpha\gamma,x\gamma\cdot z]\cdot (\alpha,x)\}=(\beta,y)\circ \{(\alpha,x)\circ [\gamma\alpha,z\alpha\cdot x]\}\Leftrightarrow (\beta,y)\circ \{\alpha\gamma\alpha,(x\gamma\cdot z)\alpha\cdot x\}=(\beta,y)\circ \{\alpha\gamma\alpha,x\gamma\alpha\cdot (z\alpha\cdot x)\}\Leftrightarrow (\beta,y)\circ \{\alpha\gamma\alpha,(x\gamma\alpha\cdot z\alpha)\cdot x\}=(\beta,y)\circ \{\alpha\gamma\alpha,x\gamma\alpha\cdot (z\alpha\cdot x)\}\Leftrightarrow \{\beta\alpha\gamma\alpha,y\alpha\gamma\alpha\cdot [(x\alpha\gamma\cdot z\alpha)\cdot x]\}=\{\beta\alpha\gamma\alpha,y\alpha\gamma\alpha\cdot [x\gamma\alpha\cdot (z\alpha\cdot x)]\}\Leftrightarrow \{y\alpha\gamma\alpha\cdot [(x\alpha\gamma\cdot z\alpha)\cdot x]\}=\{y\alpha\gamma\alpha\cdot [x\gamma\alpha\cdot (z\alpha\cdot x)]\}$.

Replace $z$ with $z\alpha^{-1}$, $y$ with $y(\alpha\gamma\alpha)^{-1}$ to get

$y\cdot [(x\alpha\gamma\cdot z)\cdot x]=y\cdot [x\alpha\gamma\cdot (z\cdot x)]$.

Put $\alpha\gamma =\delta$ to get

$y\cdot [(x\delta\cdot z)\cdot x]=y\cdot [x\delta\cdot (z\cdot x)]\Leftrightarrow zL_{x\delta}R_{x}L_{y}=zR_{x}L_{x\delta}L_{y}\Leftrightarrow L_{x\delta}R_{x}L_{y}=R_{x}L_{x\delta}L_{y}\Leftrightarrow \delta L_{x\delta}R_{x}L_{y}=\delta R_{x}L_{x\delta}L_{y}\Leftrightarrow L_{x}\delta R_{x}L_{y}=R_{x}\delta L_{x\delta}L_{y}\Leftrightarrow L_{x}R_{x}\delta L_{y}=R_{x}L_{x}\delta L_{y}\Leftrightarrow L_{x}R_{x}L_{y}\delta=R_{x}L_{x}L_{y}\delta\Leftrightarrow L_{x}R_{x}L_{y}=R_{x}L_{x}L_{y}\Leftrightarrow zL_{x}R_{x}L_{y}=zR_{x}L_{x}L_{y}$ as required.

Let $(H,\circ)$ be an $F_{40}$-algebra, and let $(\alpha,x),(\beta,y),(\gamma,z)\in H$. Then $(\beta,y)\circ \{[(\alpha,x)\circ (\alpha,x)]\circ (\gamma,z)\}=(\beta,y)\circ \{(\alpha,x)\circ [(\alpha,x)\circ (\gamma,z)]\}\Leftrightarrow (\beta,y)\circ \{[\alpha\alpha,x\alpha\cdot x]\circ (\gamma,z)\}=(\beta,y)\circ \{(\alpha,x)\circ [\alpha\gamma,x\gamma\cdot z]\}\Leftrightarrow (\beta,y)\circ [\alpha\alpha\gamma,(x\alpha\cdot x)\gamma\cdot z]=(\beta)\circ [\alpha\alpha\gamma,x\alpha\gamma\cdot (x\gamma\cdot z)]\Leftrightarrow (\beta,y)\circ [\alpha\alpha\gamma,(x\alpha\gamma\cdot x\gamma)\cdot z]=(\beta,y)\circ [\alpha\alpha\gamma,x\alpha\gamma\cdot (x\gamma\cdot z)]\Leftrightarrow \{\beta\alpha\alpha\gamma,y\alpha\alpha\gamma\cdot [(x\alpha\gamma\cdot x\gamma)\cdot z]\}=\{\beta\alpha\alpha\gamma,y\alpha\alpha\gamma\cdot [x\alpha\gamma\cdot (x\gamma\cdot z)]\}\Leftrightarrow y\alpha\alpha\gamma\cdot [(x\alpha\gamma\cdot x\alpha)\cdot z]=y\alpha\alpha\gamma\cdot [x\alpha\gamma\cdot (x\gamma\cdot z)]$. Replace $y$ with $y(\alpha\alpha\gamma)^{-1}$ to get $y\cdot [(x\alpha\gamma\cdot x\gamma)\cdot z]=y\cdot [x\alpha\gamma\cdot (x\gamma\cdot z)]$. Put $\alpha\gamma=\delta$ to get $y\cdot [(x\delta\cdot x\gamma)\cdot z]=y\cdot [x\delta\cdot (x\gamma\cdot z)]\Leftrightarrow zL_{(x\delta\cdot x\gamma)}L_{y}=zL_{x\gamma}L_{x\delta}L_{y}\Leftrightarrow L_{(x\delta\cdot x\gamma)}L_{y}=L_{x\gamma}L_{x\delta}L_{y}$. Replace $x$ with $x\gamma$ to get $L_{(x\gamma\delta\cdot x\gamma^{2})}L_{y}=L_{x\gamma^{2}}L_{x\gamma\delta}L_{y}\Leftrightarrow L_{(x\gamma\delta\cdot x)}L_{y}=L_{x}L_{x\gamma\delta}L_{y}$. Put $\gamma\delta=\theta$ to get $L_{(x\theta\cdot x)}L_{y}=L_{x}L_{x\theta}L_{y}\Leftrightarrow L_{(x\cdot x)\theta}L_{y}=L_{x}L_{x\theta}L_{y}\Leftrightarrow \theta L_{(x\cot x)\theta}L_{y}=\theta L_{x}L_{x\theta}L_{y}\Leftrightarrow \theta L_{(x\cot x)\theta}L_{y}=L_{x}\theta L_{x\theta}L_{y}\Leftrightarrow L_{x\cdot }\theta L_{y}=L_{x}L_{x}\theta L_{y}\Leftrightarrow L_{x\cdot x}L_{y}\theta=L_{x}L_{x}L_{y}\theta \Leftrightarrow L_{x\cdot x}L_{y}=L_{x}L_{x}L_{y}\Leftrightarrow zL_{x\cdot x}L_{y}=zL_{x}L_{x}L_{y}$ as required.

Let $(H,\circ)$ be an $F_{50}$-algebra, and let $(\alpha,x),(\beta,y),(\gamma,z)\in H$. Then $(\alpha,x)\circ \{(\alpha,x)\circ [(\beta,y)\circ (\gamma,z)]\}=(\alpha,x)\circ \{[(\alpha,x)\circ (\beta,y)]\circ (\gamma,z)\}\Leftrightarrow (\alpha,x)\circ \{(\alpha,x)\circ [\beta\gamma,y\gamma\cdot z]\}=(\alpha,x)\circ \{[\alpha\beta,x\beta\cdot y]\circ (\gamma,z)\}\Leftrightarrow (\alpha,x)\circ \{\alpha\beta\gamma,x\beta\gamma\cdot (y\gamma\cdot z)\}=(\alpha,x)\circ \{\alpha\beta\gamma,(x\beta\cdot y)\gamma\cdot z\}\Leftrightarrow \{\alpha\alpha\beta\gamma,x\alpha\beta\gamma\cdot (x\beta\gamma\cdot y\gamma\cdot z)\}=\{\alpha\alpha\beta\gamma,x\alpha\beta\gamma\cdot (x\beta\cdot y)\gamma\cdot z\}\Leftrightarrow \{x\alpha\beta\gamma\cdot (x\beta\gamma\cot y\gamma\cdot z)\}=\{x\alpha\beta\gamma\cdot (x\beta\gamma\cdot y\gamma)\cdot z\}$. Replace $y$ with $y\gamma^{-1}$ to get $x\alpha\beta\gamma\cdot [x\beta\gamma\cdot (y\cdot z)]=x\alpha\beta\gamma\cdot [(x\beta\gamma\cdot y)\cdot z]$. Put $\delta=\alpha\beta\gamma$ to get $x\delta\cdot [x\beta\gamma\cdot (y\cdot z)]=x\delta\cdot [(x\beta\gamma\cdot y)\cdot z]$. Put $\theta=\beta\gamma$ to get $x\delta\cdot [x\theta\cdot (y\cdot z)]=x\delta\cdot [(x\theta\cdot y)\cdot z]\Leftrightarrow yR_{z}L_{x\theta}L_{x\delta}=yL_{x\theta}R_{z}L_{x\delta}\Leftrightarrow R_{z}L_{x\theta}L_{x\delta}=L_{x\theta}R_{z}L_{x\delta}\Leftrightarrow \theta R_{z}L_{x\theta}L_{x\delta}=\theta L_{x\theta}R_{z}L_{x\delta}\Leftrightarrow R_{z}L_{x}\theta L_{x\delta}=L_{x}\theta R_{z}L_{x\delta}\Leftrightarrow R_{z}L_{x}L_{x\delta}\theta=L_{x}R_{z}L_{x\delta}\theta \Leftrightarrow \delta R_{z}L_{x}L_{x\delta}\theta=\delta L_{x}R_{z}L_{x\delta}\theta\Leftrightarrow R_{z}L_{x}\delta L_{x\delta}\theta=L_{x}R_{z}\delta L_{x\delta}\theta$

$\Leftrightarrow R_{z}L_{x}L_{x}\delta\theta=L_{x}R_{z}L_{x}\delta\theta\Leftrightarrow R_{z}L_{x}L_{x}=L_{x}R_{z}L_{x}\Leftrightarrow yR_{z}L_{x}L_{x}=yL_{x}R_{z}L_{x}$ as required.

Let $(H,\circ)$ be an $F_{53}$-algebra, and let $(\alpha,x),(\beta,y),(\gamma,z)\in H$. Then $[(\beta,y)\circ (\gamma,z)]\circ [(\alpha,x)\circ (\alpha ,x)]=(\beta,y)\circ \{[(\gamma,z)\circ (\alpha,x)]\circ (\alpha,x)\}\Leftrightarrow [\beta\gamma,y\gamma\cdot z]\circ [\alpha\alpha,x\alpha\cdot x]=(\beta,y)\circ \{[\gamma\alpha,z\alpha\cdot x]\circ (\alpha,x)\}\Leftrightarrow [\beta\gamma\alpha\alpha,(y\gamma\cdot z)\alpha\alpha\cdot (x\alpha\cdot x)]=(\beta,y)\circ [\gamma\alpha\alpha,(z\alpha\cdot x)\alpha\cdot x]\Leftrightarrow [\beta\gamma\alpha\alpha,(y\gamma\alpha\alpha\cdot z\alpha\alpha)\cdot (x\alpha\cdot x)]=\{\beta\gamma\alpha\alpha,y\gamma\alpha\alpha,\cdot [(z\alpha\alpha\cdot x\alpha)\cdot x]\}\Leftrightarrow (y\gamma\alpha\alpha\cdot z\alpha\alpha)\cdot (x\alpha\cdot x)=y\gamma\alpha\alpha\cdot [(z\alpha\alpha\cdot x\alpha)\cdot x]$. Replace $y$ with $y(\gamma\alpha\alpha)^{-1}$, $z$ with $z(\alpha\alpha)^{-1}$ to get $(y\cdot z)\cdot (x\alpha\cdot x)=y\cdot [(z\cdot x\alpha)\cdot x]\Leftrightarrow zL_{y}R_{x\alpha\cdot x}=zR_{x\alpha}R_{x}L_{y}\Leftrightarrow L_{y}R_{x\alpha\cdot x}=R_{x\alpha}R_{x}L_{y}\Leftrightarrow L_{y}R_{(x\cdot x)\alpha}=R_{x\alpha}R_{x}L_{y}\Leftrightarrow \alpha L_{y}R_{(x\cdot x)\alpha}=\alpha R_{x\alpha}R_{x}L_{y}\Leftrightarrow L_{y}\alpha R_{(x\cdot x)\alpha}=\alpha R_{x\alpha}R_{x}L_{y}\Leftrightarrow L_{y}R_{(x\cdot x)}\alpha=R_{x}\alpha R_{x}L_{y}\Leftrightarrow L_{y}R_{(x\cdot x)}\alpha=R_{x} R_{x}L_{y}\alpha\Leftrightarrow L_{y}R_{(x\cdot x)}=R_{x} R_{x}L_{y}\Leftrightarrow zL_{y}R_{(x\cdot x)}=zR_{x} R_{x}L_{y}$ as required.

Let $(H,\circ)$ be an $F_{55}$-algebra, and let $(\alpha,x),(\beta,y),(\gamma,z)\in H$. Then $\{[(\beta,y)\circ (\gamma,z)]\circ (\alpha,x)\}\circ (\alpha,x)=\{(\beta,y)\circ [(\gamma,z)\circ (\alpha,x)]\}\circ (\alpha,x)\Leftrightarrow \{[\beta\gamma,y\gamma\cdot z]\circ (\alpha,x)\}\circ (\alpha,x)=\{(\beta,y)\circ [\gamma\alpha,z\alpha\cdot x]\}\circ (\alpha,x)\Leftrightarrow \{\beta\gamma\alpha,(y\gamma\cdot z)\alpha\cdot x\}\circ (\alpha,x)=\{\beta\gamma\alpha,y\gamma\alpha\cdot (z\alpha\cdot x)\}\circ (\alpha,x)\Leftrightarrow \{\beta\gamma\alpha\alpha,[(y\gamma\alpha\cdot z\alpha)\cdot x]\alpha\cdot x\}=\{\beta\gamma\alpha\alpha,[y\gamma\alpha\cdot (z\alpha\cdot x)]\alpha\cdot x\}\Leftrightarrow \{\beta\gamma\alpha\alpha,[(y\gamma\alpha\cdot z\alpha)\cdot x]\alpha\cdot x\}=\{\beta\gamma\alpha\alpha,[y\gamma\alpha\alpha\cdot (z\alpha\alpha\cdot x\alpha)]\cdot x\}\Leftrightarrow \{\beta\gamma\alpha\alpha,[(y\gamma\alpha\alpha\cdot z\alpha\alpha)\cdot x\alpha]\cdot x\}=\{\beta\gamma\alpha\alpha,[y\gamma\alpha\alpha\cdot (z\alpha\alpha\cdot x\alpha)]\cdot x\}\Leftrightarrow [(y\gamma\alpha\alpha\cdot z\alpha\alpha)\cdot x\alpha]\cdot x=[y\gamma\alpha\alpha\cdot (z\alpha\alpha\cdot x\alpha)]\cdot x$. Replace $y$ with $y(\gamma\alpha\alpha)^{-1}$, $z$ with $z(\alpha\alpha)^{-1}$ to get $[(y\cdot z)\cdot x\alpha]\cdot x=[y\cdot (z\cdot x\alpha)]\cdot x\Leftrightarrow zL_{y}R_{x\alpha}R_{x}=zR_{x\alpha}L_{y}R_{x}\Leftrightarrow L_{y}R_{x\alpha}R_{x}=R_{x\alpha}L_{y}R_{x}\Leftrightarrow \alpha L_{y}R_{x\alpha}R_{x}=\alpha R_{x\alpha}L_{y}R_{x}\Leftrightarrow  L_{y}R_{x}\alpha R_{x}= R_{x}\alpha L_{y}R_{x}\Leftrightarrow L_{y}R_{x}R_{x}\alpha= R_{x}L_{y}R_{x}\alpha\Leftrightarrow L_{y}R_{x}R_{x}= R_{x}L_{y}R_{x} \Leftrightarrow zL_{y}R_{x}R_{x}=zR_{x}L_{y}R_{x}$ as required.

Let $(H,\circ)$ be an $F_{56}$-algebra, and let $(\alpha,x),(\beta,y),(\gamma,z)\in H$. Then $\{[(\beta,y)\circ (\gamma,z)]\circ (\alpha,x)\}\circ (\alpha,x)=(\beta,y)\circ \{[(\gamma,z)\circ (\alpha,x)]\circ (\alpha,x)\}\Leftrightarrow \{[\beta\gamma,y\gamma\cdot z]\circ (\alpha,x)\}=(\beta,y)\circ \{[\gamma\alpha,z\alpha\cdot x]\circ (\alpha,x)\}\Leftrightarrow \{\beta\gamma\alpha,(y\gamma\cdot z)\alpha\cdot x\}\circ (\alpha,x)=(\beta,y)\circ \{\gamma\alpha\alpha,(z\alpha\cdot x)\alpha\cdot x\}\Leftrightarrow \{\beta\gamma\alpha\alpha,[(y\gamma\cdot z)\alpha\cdot x]\alpha\cdot x\}=\{\beta\gamma\alpha\alpha,y\gamma\alpha\alpha\cdot [(z\alpha\cdot x)\alpha\cdot x]\}\Leftrightarrow \{\beta\gamma\alpha\alpha,[(y\gamma\alpha\cdot z\alpha)\cdot x]\alpha\cdot x\}=\{\beta\gamma\alpha\alpha,y\gamma\alpha\alpha\cdot (z\alpha\alpha\cdot x\alpha)\cdot x\}\Leftrightarrow \{\beta\gamma\alpha\alpha,[(y\gamma\alpha\cdot z\alpha)\alpha\cdot x\alpha]\cdot x\}=\{\beta\gamma\alpha\alpha,y\gamma\alpha\alpha\cdot (z\alpha\alpha\cdot x\alpha)\cdot x\}\Leftrightarrow \{\beta\gamma\alpha\alpha,[(y\gamma\alpha\alpha\cdot z\alpha\alpha)\cdot x\alpha]\cdot x\}=\{\beta\gamma\alpha\alpha,y\gamma\alpha\alpha\cdot (z\alpha\alpha\cdot x\alpha)\cdot x\}\Leftrightarrow \{[(y\gamma\alpha\alpha\cdot z\alpha\alpha)\cdot x\alpha]\cdot x\}=y\gamma\alpha\alpha\cdot \{(z\alpha\alpha\cdot x\alpha)\cdot x\}$. Replace $y$ with $y(\gamma\alpha\alpha)^{-1}$, $z$ with $z(\alpha\alpha)^{-1}$ to get $[(y\cdot z)\cdot x\alpha]\cdot x=y\cdot [(z\cdot x\alpha)\cdot x]\Leftrightarrow zL_{y}R_{x\alpha}R_{x}=zR_{x\alpha}R_{x}L_{y}\Leftrightarrow L_{y}R_{x\alpha}R_{x}=R_{x\alpha}R_{x}L_{y}\Leftrightarrow \alpha L_{y}R_{x\alpha}R_{x}=\alpha R_{x\alpha}R_{x}L_{y}\Leftrightarrow L_{y}\alpha R_{x\alpha}R_{x}=R_{x}\alpha R_{x}L_{y}\Leftrightarrow L_{y}R_{x}\alpha R_{x}=R_{x}\alpha R_{x}L_{y}\Leftrightarrow L_{y}R_{x}R_{x}\alpha=R_{x}R_{x}L_{y}\alpha\Leftrightarrow L_{y}R_{x}R_{x}=R_{x}R_{x}L_{y}\Leftrightarrow zL_{y}R_{x}R_{x}=zR_{x}R_{x}L_{y}$ as required.

Let $(H,\circ)$ be an $F_{58}$-algebra, and let $(\alpha,x),(\beta,y),(\gamma,z)\in H$. Then $\{(\beta,y)\circ [(\gamma,z)\circ (\alpha,x)]\}\circ (\alpha,x)=(\beta,y)\circ \{[(\gamma,z)\circ (\alpha,x)]\circ (\alpha,x)\}\Leftrightarrow \{(\beta,y)\circ [\gamma\alpha,z\alpha\cdot x]\}\circ (\alpha,x)=(\beta,y)\circ \{[\gamma\alpha,z\alpha\cdot x]\circ (\alpha,x)\}\Leftrightarrow \{\beta\gamma\alpha,y\gamma\alpha\cdot (z\alpha\cdot x)\}\circ (\alpha,x)=(\beta,y)\circ \{\gamma\alpha\alpha,(z\alpha\cdot x)\alpha\cdot x\}\Leftrightarrow \{\beta\gamma\alpha\alpha,[y\gamma\alpha\cdot (z\alpha\cdot x)\alpha\cdot x]\}=(\beta,y)\circ \{\gamma\alpha\alpha,(z\alpha\alpha\cdot x\alpha)\cdot x\}\Leftrightarrow \{\beta\gamma\alpha\alpha,[y\gamma\alpha\alpha\cdot (z\alpha\cdot x)\alpha]\cdot x\}=\{\beta\gamma\alpha\alpha,y\gamma\alpha\alpha\cdot (z\alpha\alpha\cdot x\alpha)\cdot x\}\Leftrightarrow \{\beta\gamma\alpha\alpha,[y\gamma\alpha\alpha\cdot (z\alpha\alpha\cdot x\alpha)]\cdot x\}=\{\beta\gamma\alpha\alpha,y\gamma\alpha\alpha\cdot (z\alpha\alpha\cdot x\alpha)\cdot x\}\Leftrightarrow [y\gamma\alpha\alpha\cdot (z\alpha\alpha\cdot x\alpha)]\cdot x=y\gamma\alpha\alpha\cdot [(z\alpha\alpha\cdot x\alpha)\cdot x]$. Replace $y$ with $y(\gamma\alpha\alpha)^{-1}$, $z$ with $z(\alpha\alpha)^{-1}$ to get $[y\cdot (z\cdot x\alpha)]\cdot x=y\cdot [(z\cdot x\alpha)\cdot x]\Leftrightarrow zR_{x\alpha}L_{y}R_{x}=zR_{x\alpha}R_{x}L_{y}\Leftrightarrow R_{x\alpha}L_{y}R_{x}=R_{x\alpha}R_{x}L_{y}\Leftrightarrow \alpha R_{x\alpha}L_{y}R_{x}=\alpha R_{x\alpha}R_{x}L_{y}\Leftrightarrow R_{x}\alpha L_{y}R_{x}=R_{x}\alpha R_{x}L_{y}\Leftrightarrow R_{x}L_{y}R_{x}\alpha=R_{x}R_{x}L_{y}\alpha\Leftrightarrow R_{x}L_{y}R_{x}=R_{x}R_{x}L_{y}\Leftrightarrow zR_{x}L_{y}R_{x}=zR_{x}R_{x}L_{y}$ as required.
\par

\noindent
{\bf Theorem 16:} Let $(G,\cdot)$ be a BCI-algebra with $A$-holomorph $(H,\circ)$ which is also a BCI-algebra such that every automorphism $\delta$ on $(G,\cdot)$ is $\lambda$-regular and $\mid\delta\mid =2$. Then $(H,\circ)$ is an $F_{i}$-algebra if and only if $(G,\cdot)$ is an $F_{i}$-algebra; where $i=4,5,6,10,20,21,25,31$.
\par

\noindent
{\bf Proof:} Let $(H,\circ)$ be an $F_{4}$-algebra, and let $(\alpha,x),(\beta,y),(\gamma,z)\in H$. Then $[(\alpha,x)\circ (\beta,y)]\circ [(\gamma,z)\circ (\alpha,x)]=(\alpha,x)\circ \{[(\beta,y)\circ (\gamma,z)]\circ (\alpha,x)\}\Leftrightarrow [(\alpha\beta,x\beta\cdot y)\circ (\gamma\alpha,z\alpha\cdot x)]=(\alpha,x)\circ [(\beta\gamma,y\gamma\cdot z)\circ (\alpha,x)]\Leftrightarrow [\alpha\beta\gamma\alpha,(x\beta\cdot y)\gamma\alpha\cdot (z\alpha\cdot x)]=(\alpha,x)\circ [\beta\gamma\alpha,(y\gamma\cdot z)\alpha\cdot x]\Leftrightarrow [\alpha\beta\gamma\alpha,(x\beta\cdot y)\gamma\alpha\cdot (z\alpha\cdot x)]=\{\alpha\beta\gamma\alpha,x\beta\gamma\alpha\cdot [(y\gamma\cdot z)\alpha\cdot x]\}\Leftrightarrow [\alpha\beta\gamma\alpha,(x\beta\gamma\alpha\cdot y\gamma\alpha)\cdot (z\alpha\cdot x)]=\{\alpha\beta\gamma\alpha,x\beta\gamma\alpha\cdot [(y\gamma\alpha\cdot z\alpha)\cdot x]\}\Leftrightarrow (x\beta\gamma\alpha\cdot y\gamma\alpha)\cdot (z\alpha\cdot x)x\beta\gamma\alpha\cdot [(y\gamma\alpha\cdot z\alpha)\cdot x]$. Replace $z$ with $z\alpha^{-1}$, $y$ with $y(\gamma\alpha)^{-1}$ to get $(x\beta\gamma\alpha\cdot y)\cdot (z\cdot x)=x\beta\gamma\gamma\alpha\cdot [(y\cdot z)\cdot x]\Leftrightarrow (x\delta\cdot y)\cdot (z\cdot x)=x\delta\cdot [(y\cdot z)\cdot x]$(where $\delta=\beta\gamma\alpha$) $\Leftrightarrow yL_{x\delta}R_{z\cdot x}=yR_{z}R_{x}L_{x\delta}\Leftrightarrow L_{x\delta}R_{z\cdot x}=R_{z}R_{x}L_{x\delta}\Leftrightarrow L_{x}\delta R_{z\cdot x}=R_{z}R_{x}L_{x}\delta\Leftrightarrow L_{x}R_{z\cdot x}\delta=R_{z}R_{x}L_{x}\delta\Leftrightarrow L_{x}R_{z\cdot x}=R_{z}R_{x}L_{x}\Leftrightarrow yL_{x}R_{z\cdot x}=yR_{z}R_{x}L_{x}$ as required.

Let $(H,\circ)$ be an $F_{5}$-algebra, and let $(\alpha,x),(\beta,y),(\gamma,z)\in H$. Then $\{[(\alpha,x)\circ (\beta,y)]\circ (\gamma,z)\}\circ (\alpha,x)=\{(\alpha,x)\circ [(\beta,y)\circ (\gamma,z)]\}\circ (\alpha,x)\Leftrightarrow \{[\alpha\beta,x\beta\cdot y]\circ (\gamma,z)\}\circ (\alpha,x)=\{(\alpha,x)\circ [\beta\gamma,y\gamma\cdot z]\}\circ (\alpha,x)\Leftrightarrow \{\alpha\beta\gamma,(x\beta\cdot y)\gamma\cdot z\}\circ (\alpha,x)=\{\alpha\beta\gamma,(x\beta\gamma)\cdot (y\gamma\cdot z)\}\circ (\alpha,x)\Leftrightarrow \{\alpha\beta\gamma\alpha,[(x\beta\cdot y)\gamma\cdot z]\alpha\cdot x\}=\{\alpha\beta\gamma\alpha,[(x\beta\gamma)\cdot (y\gamma\cdot z)]\alpha\cdot x\}\Leftrightarrow \{\alpha\beta\gamma\alpha,[(x\beta\gamma\cdot y\gamma)\cdot z]\alpha\cdot x\}=\{\alpha\beta\gamma\alpha,[(x\beta\gamma\alpha)\cdot (y\gamma\cdot z)\alpha]\cdot x\}\Leftrightarrow \{\alpha\beta\gamma\alpha,[(x\beta\gamma\cdot y\gamma)\cdot z]\alpha\cdot x\}=\{\alpha\beta\gamma\alpha,[(x\beta\gamma\alpha)\cdot (y\gamma\alpha\cdot z\alpha)]\cdot x\}\Leftrightarrow \{\alpha\beta\gamma\alpha,[(x\beta\gamma\cdot y\gamma)\alpha\cdot z\alpha]\cdot x\}=\{\alpha\beta\gamma\alpha,[(x\beta\gamma\alpha)\cdot (y\gamma\alpha\cdot z\alpha)]\cdot x\}\Leftrightarrow \{\alpha\beta\gamma\alpha,[(x\beta\gamma\alpha\cdot y\gamma\alpha)\cdot z\alpha]\cdot x\}=\{\alpha\beta\gamma\alpha,[(x\beta\gamma\alpha)\cdot (y\gamma\alpha\cdot z\alpha)]\cdot x\}\Leftrightarrow [(x\beta\gamma\alpha\cdot y\gamma\alpha)\cdot z\alpha]\cdot x=[(x\beta\gamma\alpha)\cdot (y\gamma\alpha\cdot z\alpha)]\cdot x$.

Replace $z$ with $z\alpha^{-1}$, $y$ with $y(\gamma\alpha)^{-1}$ to get

$[(x\beta\gamma\alpha\cdot y)\cdot z]\cdot x=[(x\beta\gamma\alpha)\cdot (y\cdot z)]\cdot x$.

Put $\delta=\beta\gamma\alpha$ to get

$[(x\delta\cdot y)\cdot z]\cdot x=[x\delta\cdot (y\cdot z)]\cdot x\Leftrightarrow yL_{x\delta}R_{z}R_{x}=yR_{z}L_{x\delta}R_{x}\Leftrightarrow L_{x\delta}R_{z}R_{x}=R_{z}L_{x\delta}R_{x}\Leftrightarrow \delta L_{x\delta}R_{z}R_{x}=\delta R_{z}L_{x\delta}R_{x}\Leftrightarrow L_{x}R_{z}R_{x}\delta=R_{z}L_{x}R_{z}\delta\Leftrightarrow L_{x}R_{z}R_{x}=R_{z}L_{x}R_{z}\Leftrightarrow yL_{x}R_{z}R_{x}=yR_{z}L_{x}R_{z}$ as required.

Let $(H,\circ)$ be an $F_{6}$-algebra, and let $(\alpha,x),(\beta,y),(\gamma,z)\in H$. Then $\{[(\alpha,x)\circ (\beta,y)]\circ (\gamma,z)\}\circ (\alpha,x)=(\alpha,x)\circ \{(\beta,y)\circ [(\gamma,z)\circ (\alpha,x)]\}\Leftrightarrow \{[\alpha\beta,x\beta\cdot y]\circ (\gamma,z)\}\circ (\alpha,x)=(\alpha,x)\circ \{(\beta,y)\circ [\gamma\alpha,z\alpha\cdot x]\}\Leftrightarrow [\alpha\beta\gamma,(x\beta\cdot y)\gamma\cdot z]\circ (\alpha,x)=(\alpha,x)\circ [\beta\gamma\alpha,y\gamma\alpha\cdot (z\alpha\cdot x)]\Leftrightarrow \{\alpha\beta\gamma\alpha,[(x\beta\cdot y)\gamma\cdot z]\alpha\cdot x\}=\{\alpha\beta\gamma\alpha,x\beta\gamma\alpha\cdot [y\gamma\alpha\cdot (z\alpha\cdot x)]\}\Leftrightarrow \{\alpha\beta\gamma\alpha,[(x\beta\gamma\cdot y\gamma)\cdot z]\alpha\cdot x\}=\{\alpha\beta\gamma\alpha,x\beta\gamma\alpha\cdot [y\gamma\alpha\cdot (z\alpha\cdot x)]\}\Leftrightarrow \{\alpha\beta\gamma\alpha,[(x\beta\gamma\cdot y\gamma)\alpha\cdot z\alpha]\cdot x\}=\{\alpha\beta\gamma\alpha,x\beta\gamma\alpha\cdot [y\gamma\alpha\cdot (z\alpha\cdot x)]\}\Leftrightarrow \{\alpha\beta\gamma\alpha,[(x\beta\gamma\alpha\cdot y\gamma\alpha)\cdot z\alpha]\cdot x\}=\{\alpha\beta\gamma\alpha,x\beta\gamma\alpha\cdot [y\gamma\alpha\cdot (z\alpha\cdot x)]\}\Leftrightarrow [(x\beta\gamma\alpha\cdot y\gamma\alpha)\cdot z\alpha]\cdot x=\{x\beta\gamma\alpha\cdot [y\gamma\alpha\cdot (z\alpha\cdot x)]\}$. Replace $z$ with $z\alpha^{-1}$, $y$ with $y(\gamma\alpha)^{-1}$ to get $[(x\beta\gamma\alpha\cdot y)\cdot z]\cdot x=\{x\beta\gamma\alpha\cdot [y\cdot (z\cdot x)]\}$. Put $\delta=\beta\gamma\alpha$ to get $[(x\delta\cdot y)\cdot z]\cdot x=x\delta\cdot [y\cdot (z\cdot x)]\Leftrightarrow yL_{x\delta}R_{z}R_{x}=yR_{z\cdot x}L_{x\delta}\Leftrightarrow L_{x\delta}R_{z}R_{x}=R_{z\cdot x}L_{x\delta}\Leftrightarrow \delta L_{x\delta}R_{z}R_{x}=\delta R_{z\cdot x}L_{x\delta}\Leftrightarrow L_{x}\delta R_{z}R_{x}=\delta R_{z\cdot x}L_{x\delta}\Leftrightarrow L_{x}R_{z}R_{x}\delta=R_{z\cdot x}\delta L_{x\delta}\Leftrightarrow L_{x}R_{z}R_{x}\delta=R_{z\cdot x}L_{x}\delta\Leftrightarrow L_{x}R_{z}R_{x}=R_{z\cdot x}L_{x}\Leftrightarrow yL_{x}R_{z}R_{x}=yR_{z\cdot x}L_{x}$ as required.

Let $(H,\circ)$ be an $F_{10}$-algebra, and let $(\alpha,x),(\beta,y),(\gamma,z)\in H$. Then $(\alpha,x)\circ \{(\beta,y)\circ [(\gamma,z)\circ (\alpha,x)]\}=(\alpha,x)\circ \{[(\beta,y)\circ (\gamma,z)]\circ (\alpha,x)\}\Leftrightarrow (\alpha,x)\circ \{(\beta,y)\circ [\gamma\alpha,z\alpha\cdot x]\}=(\alpha,x)\circ \{[\beta\gamma,y\gamma\cdot z]\circ (\alpha,x)\}\Leftrightarrow (\alpha,x)\circ \{\beta\gamma\alpha,y\gamma\alpha\cdot (z\alpha\cdot x)\}=(\alpha,x)\circ \{\beta\gamma\alpha,(y\gamma\cdot z)\alpha\cdot x\}\Leftrightarrow \{\alpha\beta\gamma\alpha,x\beta\gamma\alpha\cdot [y\gamma\alpha\cdot (z\alpha\cdot x)]\}=\{\alpha\beta\gamma\alpha,x\beta\gamma\alpha\cdot (y\gamma\cdot z)\alpha\cdot x\}\Leftrightarrow \{\alpha\beta\gamma\alpha,x\beta\gamma\alpha\cdot [y\gamma\alpha\cdot (z\alpha\cdot x)]\}=\{\alpha\beta\gamma\alpha,x\beta\gamma\alpha\cdot (y\gamma\alpha\cdot z\alpha)\cdot x\}\Leftrightarrow \{x\beta\gamma\alpha\cdot [y\gamma\alpha\cdot (z\alpha\cdot x)]\}=\{x\beta\gamma\alpha\cdot (y\gamma\alpha\cdot z\alpha)\cdot x\}$. Replace $z$ with $z\alpha^{-1}$, $y$ with $y(\gamma\alpha)^{-1}$ to get $\{x\beta\gamma\alpha\cdot [y\cdot (z\cdot x)]\}=\{x\beta\gamma\alpha\cdot (y\cdot z)\cdot x\}$. Put $\delta=\beta\gamma\alpha$ to get $x\delta\cdot [y\cdot (z\cdot x)]=x\delta\cdot [(y\cdot z)\cdot x]\Leftrightarrow yR_{z\cdot x}L_{x\delta}=yR_{z}R_{x}L_{x\delta}\Leftrightarrow R_{z\cdot x}L_{x\delta}=R_{z}R_{x}L_{x\delta}\Leftrightarrow \delta R_{z\cdot x}L_{x\delta}=\delta R_{z}R_{x}L_{x\delta}\Leftrightarrow R_{z\cdot x}\delta L_{x\delta}=R_{z}R_{x}\delta L_{x\delta}\Leftrightarrow R_{z\cdot x}L_{x}\delta=R_{z}R_{x} L_{x}\delta\Leftrightarrow R_{z\cdot x}L_{x}=R_{z}R_{x} L_{x}\Leftrightarrow yR_{z\cdot x}L_{x}=yR_{z}R_{x} L_{x}$ as required.

Let $(H,\circ)$ be an $F_{20}$-algebra, and let $(\alpha,x),(\beta,y),(\gamma,z)\in H$. Then $(\alpha,x)\circ \{[(\beta,y)\circ (\alpha,x)]\circ (\gamma,z)\}=(\alpha,x)\circ \{(\beta,y)\circ [(\alpha,x)\circ (\gamma,z)]\}\Leftrightarrow (\alpha,x)\circ \{[\beta\alpha,y\alpha\cdot x]\circ (\gamma,z)\}=(\alpha,x)\circ \{(\beta,y)\circ [\alpha\gamma,x\gamma\cdot z]\}\Leftrightarrow (\alpha,x)\circ \{\beta\alpha\gamma,(y\alpha\cdot x)\gamma\cdot z\}=(\alpha,x)\circ \{\beta\alpha\gamma,y\alpha\gamma\cdot (x\gamma\cdot z)\}\Leftrightarrow \{\alpha\beta\alpha\gamma,x\beta\alpha\gamma\cdot [(y\alpha\cdot x)\gamma\cdot z]\}=\{\alpha\beta\alpha\gamma,x\beta\alpha\gamma\cdot [y\alpha\gamma\cdot (x\gamma\cdot z)]\}\Leftrightarrow \{\alpha\beta\alpha\gamma,x\beta\alpha\gamma\cdot (y\alpha\gamma\cdot x\gamma)\cdot z\}=\{\alpha\beta\alpha\gamma,x\beta\alpha\gamma\cdot [y\alpha\gamma\cdot (x\gamma\cdot z)]\}\Leftrightarrow \{x\beta\alpha\gamma\cdot (y\alpha\gamma\cdot x\gamma)\cdot z\}=\{x\beta\alpha\gamma\cdot [y\alpha\gamma\cdot (x\gamma\cdot z)]\}$. Replace $y$ with $y(\alpha\gamma)^{-1}$ to get $\{x\beta\alpha\gamma\cdot (y\cdot x)\cdot z\}=\{x\beta\alpha\gamma\cdot [y\cdot (x\cdot z)]\}$. Put $\delta=\beta\alpha\gamma$ to get $x\delta\cdot [(y\cdot x)\cdot z]=\{x\delta\cdot [y\cdot (x\cdot z)]\}\Leftrightarrow yR_{x}R_{z}L_{x\delta}=yR_{x\cdot z}L_{x\delta}\Leftrightarrow R_{x}R_{z}L_{x\delta}=R_{x\cdot z}L_{x\delta}\Leftrightarrow \delta R_{x}R_{z}L_{x\delta}=\delta R_{x\cdot z}L_{x\delta}\Leftrightarrow R_{x}R_{z}\delta L_{x\delta}=\delta R_{x\cdot z}\delta L_{x\delta}\Leftrightarrow R_{x}R_{z}L_{x}\delta =R_{x\cdot z} L_{x}\delta\Leftrightarrow R_{x}R_{z}L_{x}=R_{x\cdot z} L_{x}\Leftrightarrow yR_{x}R_{z}L_{x}=yR_{x\cdot z} L_{x}$ as required.

Let $(H,\circ)$ be an $F_{21}$-algebra, and let $(\alpha,x),(\beta,y),(\gamma,z)\in H$. Then $[(\beta,y)\circ (\alpha,x)]\circ [(\gamma,z)\circ (\alpha,x)]=\{[(\beta,y)\circ (\alpha,x)]\circ (\gamma,z)\}\circ (\alpha,x)\Leftrightarrow (\beta\alpha,y\alpha\cdot x)\circ (\gamma\alpha,z\alpha\cdot x)=\{[\beta\alpha,y\alpha\cdot x]\circ (\gamma,z)\}\circ (\alpha,x)\Leftrightarrow [\beta\alpha\gamma\alpha,(y\alpha\cdot x)\gamma\alpha\cdot (z\alpha\cdot x)]=[\beta\alpha\gamma,(y\alpha\cdot x)\gamma\cdot z]\circ (\alpha,x)\Leftrightarrow [\beta\alpha\gamma\alpha,(y\alpha\cdot x)\gamma\alpha\cdot (z\alpha\cdot x)]=\{\beta\alpha\gamma\alpha,[(y\alpha\cdot x)\gamma\cdot z]\alpha\cdot x\}\Leftrightarrow [\beta\alpha\gamma\alpha,(y\alpha\gamma\alpha\cdot x\gamma\alpha)\cdot (z\alpha\cdot x)]=\{\beta\alpha\gamma\alpha,[(y\alpha\gamma\cdot x\gamma)\cdot z]\alpha\cdot x\}\Leftrightarrow [\beta\alpha\gamma\alpha,(y\alpha\gamma\alpha\cdot x\gamma\alpha)\cdot (z\alpha\cdot x)]=\{\beta\alpha\gamma\alpha,[(y\alpha\gamma\cdot x\gamma)\alpha\cdot z\alpha]\cdot x\}\Leftrightarrow [\beta\alpha\gamma\alpha,(y\alpha\gamma\alpha\cdot x\gamma\alpha)\cdot (z\alpha\cdot x)]=\{\beta\alpha\gamma\alpha,[(y\alpha\gamma\alpha\cdot x\gamma\alpha)\cdot z\alpha]\cdot x\}\Leftrightarrow (y\alpha\gamma\alpha\cdot x\gamma\alpha)\cdot (z\alpha\cdot x)=[(y\alpha\gamma\alpha\cdot x\gamma\alpha)\cdot z\alpha]\cdot x$. Replace $y$ with $y(\alpha\gamma\alpha)^{-1}$, $z$ with $z\alpha^{-1}$ to get $(y\cdot x\gamma\alpha)\cdot (z\cdot x)=[(y\cdot x\gamma\alpha)\cdot z]\cdot x$. Put $\delta=\gamma\alpha$ to get $(y\cdot x\delta)\cdot (z\cdot x)=[(y\cdot x\delta)\cdot z]\cdot x\Leftrightarrow yR_{x\delta}R_{z\cdot x}=yR_{x\delta}R_{z}R_{x}\Leftrightarrow R_{x\delta}R_{z\cdot x}=R_{x\delta}R_{z}R_{x}\Leftrightarrow \delta R_{x\delta}R_{z\cdot x}=\delta R_{x\delta}R_{z}R_{x}\Leftrightarrow R_{x}\delta R_{z\cdot x}= R_{x}\delta R_{z}R_{x}\Leftrightarrow R_{x}R_{z\cdot x}\delta=R_{x}R_{z}R_{x}\delta\Leftrightarrow R_{x}R_{z\cdot x}=R_{x}R_{z}R_{x}\Leftrightarrow yR_{x}R_{z\cdot x}=yR_{x}R_{z}R_{x}$ as required.

Let $(H,\circ)$ be an $F_{25}$-algebra, and let $(\alpha,x),(\beta,y),(\gamma,z)\in H$. Then $\{[(\beta,y)\circ (\alpha,x)]\circ (\gamma,z)\}\circ (\alpha,x)=\{(\beta,y)\circ [(\alpha,x)\circ (\gamma,z)]\}\circ (\alpha,x)\Leftrightarrow \{[\beta\alpha,y\alpha\cdot x]\circ (\gamma,z)\}\circ (\alpha,x)=\{(\beta,y)\circ [\alpha\gamma,x\gamma\cdot z]\}\circ (\alpha,x)\Leftrightarrow \{[\beta\alpha\gamma,(y\alpha\cdot x)\gamma\cdot z]\}\circ (\alpha,x)=[\beta\alpha\gamma\cdot (x\gamma\cdot z)]\circ (\alpha,x)\Leftrightarrow \{[\beta\alpha\gamma,(y\alpha\gamma\cdot x\gamma)\cdot z]\}\circ (\alpha,x)=[\beta\alpha\gamma,y\alpha\gamma\cdot (x\gamma\cdot z)]\circ (\alpha,x)\Leftrightarrow \{\beta\alpha\gamma\alpha,[(y\alpha\gamma\cdot x\gamma)\cdot z]\alpha\cdot x\}=\{\beta\alpha\gamma\alpha,[y\alpha\gamma\cdot (x\gamma\cdot z)]\alpha\cdot x\}\Leftrightarrow \{\beta\alpha\gamma\alpha,[(y\alpha\gamma\cdot x\gamma)\alpha\cdot z\alpha]\cdot x\}=\{\beta\alpha\gamma\alpha,[y\alpha\gamma\alpha\cdot (x\gamma\cdot z)\alpha]\cdot x\}\Leftrightarrow \{\beta\alpha\gamma\alpha,[(y\alpha\gamma\alpha\cdot x\gamma\alpha)\cdot z\alpha]\cdot x\}=\{\beta\alpha\gamma\alpha,[y\alpha\gamma\alpha\cdot (x\gamma\alpha\cdot z\alpha)]\cdot x\}\Leftrightarrow [(y\alpha\gamma\alpha\cdot x\gamma\alpha)\cdot z\alpha]\cdot x=[y\alpha\gamma\alpha\cdot (x\gamma\alpha\cdot z\alpha)]\cdot x$. Replace $y$ with $y(\alpha\gamma\alpha)^{-1}$, $z$ with $z\alpha^{-1}$ to get $[(y\cdot x\gamma\alpha)\cdot z]\cdot x=[y\cdot (x\gamma\alpha\cdot z)]\cdot x$. Put $\delta=\gamma\alpha$ to get $[(y\cdot x\delta)\cdot z]\cdot x=[y\cdot (x\delta\cdot z)]\cdot x$. Replace $x$ with $x\delta$ to get $[(y\cdot x\delta\delta)\cdot z]\cdot x\delta=[y\cdot (x\delta\delta\cdot z)]\cdot x\delta\Leftrightarrow [(y\cdot x)\cdot z]\cdot x\delta=[y\cdot (x\cdot z)]\cdot x\delta\Leftrightarrow yR_{x}R_{z}R_{x\delta}=yR_{x\cdot z}R_{x\delta}\Leftrightarrow R_{x}R_{z}R_{x\delta}=R_{x\cdot z}R_{x\delta}\Leftrightarrow \delta R_{x}R_{z}R_{x\delta}=\delta R_{x\cdot z}R_{x\delta}\Leftrightarrow R_{x}R_{z}\delta R_{x\delta}=R_{x\cdot z}\delta R_{x\delta}\Leftrightarrow R_{x}R_{z}R_{x}\delta=R_{x\cdot z}R_{x}\delta\Leftrightarrow R_{x}R_{z}R_{x}=R_{x\cdot z}R_{x}\Leftrightarrow yR_{x}R_{z}R_{x}=yR_{x\cdot z}R_{x}$ as required.

Let $(H,\circ)$ be an $F_{31}$-algebra, and let $(\alpha,x),(\beta,y),(\gamma,z)\in H$. Then $[(\beta,y)\circ (\alpha,x)]\circ [(\alpha,x)\circ (\gamma,z)]=\{[(\beta,y)\circ (\alpha,x)]\circ (\alpha,x)\}\circ (\gamma,z)\Leftrightarrow [\beta\alpha,y\alpha\cdot x]\circ [\alpha\gamma,x\gamma\cdot z]=\{[\beta\alpha,y\alpha\cdot x]\circ (\alpha,x)\}\circ (\gamma,z)\Leftrightarrow [\beta\alpha\alpha\gamma,(y\alpha\cdot x)\alpha\gamma\cdot (x\gamma\cdot z)]=\{[\beta\alpha\alpha,(y\alpha\cdot x)\alpha\cdot x]\}\circ (\gamma,z)\Leftrightarrow [\beta\alpha\alpha\gamma,(y\alpha\alpha\gamma\cdot x\alpha\gamma)\cdot (x\gamma\cdot z)]=\{[\beta\alpha\alpha,(y\alpha\alpha\cdot x\alpha)\cdot x]\}\circ (\gamma,z)\Leftrightarrow [\beta\alpha\alpha\gamma,(y\alpha\alpha\gamma\cdot x\alpha\gamma)\cdot (x\gamma\cdot z)]=\{\beta\alpha\alpha\gamma,[(y\alpha\alpha\cdot x\alpha)\cdot x]\gamma cdot z\}\Leftrightarrow [\beta\alpha\alpha\gamma,(y\alpha\alpha\gamma\cdot x\alpha\gamma)\cdot (x\gamma\cdot z)]=\{\beta\alpha\alpha\gamma,[(y\alpha\alpha\gamma\cdot x\alpha\gamma)\cdot x\gamma\cdot z]\}\Leftrightarrow (y\alpha\alpha\gamma\cdot x\alpha\gamma)\cdot (x\gamma\cdot z)=[(y\alpha\alpha\gamma\cdot x\alpha\gamma)\cdot x\gamma]\cdot z$. Replace $y$ with $y(\alpha\alpha)^{-1}$ to get $(y\cdot x\alpha\gamma)\cdot (x\gamma\cdot z)=(y\cdot x\alpha\gamma)\cdot x\gamma\cdot z$. Put $\delta=\alpha\gamma$ to get $(y\cdot x\delta)\cdot (x\gamma\cdot z)=[(y\cdot x\delta)\cdot x\gamma]\cdot z\Leftrightarrow yR_{x\delta}R_{x\gamma\cdot z}=yR_{x\delta}R_{x\gamma}R_{z}\Leftrightarrow R_{x\delta}R_{x\gamma\cdot z}=R_{x\delta}R_{x\gamma}R_{z}$. Put $x=x\delta$ to get $R_{x\delta^{2}}R_{x\delta\gamma\cdot z}=R_{x\delta^{2}}R_{x\delta\gamma}R_{z}\Leftrightarrow R_{x}R_{x\delta\gamma\cdot z}=R_{x}R_{x\delta\gamma}R_{z}$. Put $\delta\gamma=\theta$ to get $R_{x}R_{x\theta\cdot z}=R_{x}R_{x\theta}R_{z}\Leftrightarrow R_{x}R_{(x\cdot z)\theta}=R_{x}R_{x\theta}R_{z}\Leftrightarrow \theta R_{x}R_{(x\cdot z)\theta}=\theta R_{x}R_{x\theta}R_{z}\Leftrightarrow R_{x}\theta R_{(x\cdot z)\theta}=R_{x}\theta R_{x\theta}R_{z}\Leftrightarrow R_{x}R_{x\cdot z}\theta=R_{x}R_{x}\theta R_{z}\Leftrightarrow R_{x}R_{x\cdot z}\theta=R_{x}R_{x}R_{z}\theta\Leftrightarrow R_{x}R_{x\cdot z}=R_{x}R_{x}R_{z}\Leftrightarrow yR_{x}R_{x\cdot z}=yR_{x}R_{x}R_{z}$ as required.
\par

\noindent
{\bf Theorem 17:} Let $(G,\cdot)$ be a BCI-algebra with $A$-holomorph $(H,\circ)$ which is also a BCI-algebra. Then $(H,\circ)$ is an $F_{i}$-algebra if and only if $(G,\cdot)$ is an $F_{i}$-algebra; where $i=42,54$.
\par

\noindent
{\bf Proof:} Let $(H,\circ)$ be an $F_{42}$-algebra, and let $(\alpha,x),(\beta,y),(\gamma,z)\in H$. Then $(I,0)\circ [(\beta,y)\circ (\gamma,z)]=[(I,0)\circ (\beta,y)]\circ (\gamma,z)\Leftrightarrow (I,0)\circ [\beta\gamma,y\gamma\cdot z]=[I\beta,0\beta\cdot y]\circ (\gamma,z)\Leftrightarrow [I\beta\gamma,0\beta\gamma\cdot (y\gamma\cdot z)]=[I\beta\gamma,(0\beta\cdot y)\gamma\cdot z]\Leftrightarrow [\beta\gamma,0\beta\gamma\cdot (y\gamma\cdot z)]=[\beta\gamma,(0\beta\gamma\cdot y\gamma)\cdot z]\Leftrightarrow 0\beta\gamma\cdot (y\gamma\cdot z)=(0\beta\gamma\cdot y\gamma)\cdot z\Leftrightarrow 0\cdot (y\gamma\cdot z)=(0\cdot y\gamma)\cdot z\Leftrightarrow y\gamma R_{z}L_{0}=y\gamma L_{0}R_{z}\Leftrightarrow \gamma R_{z}L_{0}=\gamma L_{0}R_{z}\Leftrightarrow R_{z}L_{0}=L_{0}R_{z}\Leftrightarrow yR_{z}L_{0}=yL_{0}R_{z}$ as required.

Let $(H,\circ)$ be an $F_{54}$-algebra, and let $(\beta,y),(\gamma,z)\in H$. Then $(\beta,y)\circ (\gamma,z)=(\beta,y)\circ (\gamma,z)\Leftrightarrow (\beta\gamma,y\gamma\cdot z)=(\beta\gamma,y\gamma\cdot z)$. Put $\beta\gamma=\delta$ to get $(\delta,y\gamma\cdot z)=(\delta,y\gamma\cdot z)\Leftrightarrow y\gamma\cdot z=y\gamma\cdot z\Leftrightarrow zL_{y\gamma}=zL_{y\gamma}\Leftrightarrow L_{y\gamma}=L_{y\gamma}\Leftrightarrow \gamma L_{y\gamma}=\gamma L_{y\gamma}\Leftrightarrow L_{y}\gamma=L_{y}\gamma\Leftrightarrow L_{y}=L_{y}\Leftrightarrow zL_{y}=zL_{y}$ as required.
\par




\vskip 0.3 cm
{\small{
\centerline{REFERENCES}
}}
\scriptsize{

}
\par
\noindent 
\tiny{DEPARTMENT OF MATHEMATICS,
FEDERAL UNIVERSITY OF AGRICULTURE, ABEOKUTA, NIGERIA
\par
\noindent 
{\it E-mail address}: {\tt emmanuelilojide@gmail.com, ilojidee@unaab.edu.ng}
\par
\noindent 
DEPARTMENT OF MATHEMATICS,
OBAFEMI AWOLOWO UNIVERSITY, ILE-IFE, NIGERIA
\par
\noindent 
{\it E-mail addresses}: {\tt jaiyeolatemitope@yahoo.com, tjayeola@oauife.edu.ng}

\noindent
DEPARTMENT OF MATHEMATICS,
OBAFEMI AWOLOWO UNIVERSITY, ILE-IFE, NIGERIA
\par
\noindent
{\it E-mail addresses}: {\tt memudu.olatinwo@gmail.com}
\par
\end{document}